%%%%%%%%%%%%%%%%AMS TeXfile%%%%%%%%%%%%%%%%%%%%%%%%%%%%%%%%%%%%%%%%%%%%
\input amstex.tex
\documentstyle{amsppt}
\magnification=1200
\baselineskip=13pt
\hsize=6.5truein
\vsize=8.9truein
%%%%%%%%%%%%%%%%% M a c r o s %%%%%%%%%%%%%%%%%%%%%%%%%%%%%%%%%%%%%%%%%%%
%Section, Theorem and Formula Numbering%%%%%%%%%%%%%%%%%%%%%%%%%%%%%%%%%%
%%%%%%%%%%%%%%%%%%%%%%%%%%%%%%%%%%%%%%%%%%%%%%%%%%%%%%%%%%%%%%%%%%%%%%%%%
\countdef\sectionno=1
\countdef\eqnumber=10
\countdef\theoremno=11
\countdef\countrefno=12
\countdef\cntsubsecno=13
\sectionno=0
\def\newsection{\global\advance\sectionno by 1
                \global\eqnumber=1
                \global\theoremno=1
                \global\cntsubsecno=0
                \the\sectionno.\ }

\def\newsubsection#1{\global\advance\cntsubsecno by 1
                     \xdef#1{{\S\the\sectionno.\the\cntsubsecno}}
                     \ \the\sectionno.\the\cntsubsecno.}

\def\theoremname#1{\the\sectionno.\the\theoremno
                   \xdef#1{{\the\sectionno.\the\theoremno}}
                   \global\advance\theoremno by 1}

\def\eqname#1{\the\sectionno.\the\eqnumber
              \xdef#1{{\the\sectionno.\the\eqnumber}}
              \global\advance\eqnumber by 1}

\def\thmref#1{#1}

\global\countrefno=1

\def\refno#1{\xdef#1{{\the\countrefno}}\global\advance\countrefno by 1}

%%%%%%%%%%%%%%%%%%%%%%%%%%%%Abbreviations%%%%%%%%%%%%%%%%%%%%%%%%%%%%%%%%
\def\R{{\Bbb R}}
\def\N{{\Bbb N}}
\def\C{{\Bbb C}}

\def\Zp{{\Bbb Z}_+}
\def\hf{{1\over 2}}
\def\hZp{{1\over 2}\Zp}
\def\A{{\Cal A}_q(SU(2))}
\def\Hi{\ell^2(\Zp)}
\def\a{\alpha}
\def\b{\beta}
\def\g{\gamma}
\def\d{\delta}
\def\s{\sigma}
\def\t{\tau}
\def\l{\lambda}
\def\m{\mu}

\def\r{\rho_{\t ,\s}}

\def\rti{\rho_{\t ,\infty}}
\def\p{\pi}
\def\vp{\varphi}
%%%%%%%%%%%%%%%%%%%%%%%%%Reference Numbering%%%%%%%%%%%%%%%%%%%%%%%%%
\refno{\Aske}
\refno{\AskeIqH}
\refno{\AskeI}
\refno{\AskeRS}
\refno{\AskeW}
\refno{\Bere}
\refno{\Bres}
\refno{\CharP}
\refno{\DijkN}
\refno{\Domb}
\refno{\Dunk}
\refno{\ErdeHTF}
\refno{\FlorLV}
\refno{\FlorVJMP}
\refno{\FlorVPLA}
\refno{\FlorVLMP}
\refno{\FlorVJGTP}
\refno{\FlorVAP}
\refno{\FlorVCJP}
\refno{\FlorVJCAM}
\refno{\FlorVEsterel}
\refno{\Flor}
\refno{\FlorK}
\refno{\GaspR}
\refno{\GrozK}
\refno{\KalnMM}
\refno{\KalnM}
\refno{\KalnMiMuJMP}
\refno{\KalnMiMuSIAM}
\refno{\KoelITSF}
\refno{\KoelSIAM}
\refno{\KoelDMJ}
\refno{\KoelCJM}
\refno{\KoelIM}
\refno{\KoelAAM}
\refno{\KoelJCAM}
\refno{\KoelCJMs}
\refno{\KoelS}
\refno{\KoorIM}
\refno{\KoorOPTA}
\refno{\KoorSIAM}
\refno{\KoorZSE}
\refno{\KoorS}
\refno{\Mill}
\refno{\MasuMNNU}
\refno{\Noum}
\refno{\NoumMPJA}
\refno{\NoumMCMP}
\refno{\NoumMDMJ}
\refno{\NoumMCM}
\refno{\NoumMLNM}
\refno{\RahmCJM}
\refno{\RahmPAMS}
\refno{\RahmV}
\refno{\StanSIAM}
\refno{\StanGD}
\refno{\StanAKS}
\refno{\Swar}
\refno{\VaksS}
\refno{\VAsscK}
\refno{\Vile}
\refno{\VileK}
\refno{\Wats}
\refno{\Woro}

%%%%%Beginning of the text%%%%%%%%%%%%%%%%%%%%%%%%%%%%%%%%%%%%%%%%%%%%
\topmatter
\title Addition formulas for $q$-special functions\endtitle
\author H.T. Koelink\endauthor
\address Vakgroep Wiskunde, Universiteit van Amsterdam,
Plantage Muidergracht 24, 1018 TV Amsterdam, the Netherlands\endaddress
\email koelink\@fwi.uva.nl\endemail
%\date Version of \the\day-\the\month-\the\year\enddate
\thanks Research supported by a Fellowship of the Research Council of the
Katholieke Universiteit Leuven and the Netherlands Organization for Scientific
Research (NWO) under project number 610.06.100.\endthanks
\subjclass 33D80, 33D45, 33D55, 17B37\endsubjclass
\abstract  A general addition formula for a two-parameter
family of Askey-Wilson polynomials is derived from the quantum $SU(2)$
group theoretic interpretation. This formula contains most of the
previously known addition formulas for $q$-Legendre polynomials as special
or limiting cases. A survey of the literature on addition formulas for
$q$-special functions using quantum groups and quantum algebras is given.
\endabstract
\endtopmatter
\document

%%%%%%%%%%%%%%%%%%%%%%%%%%%%%%%%%%%%%%%%%%%%%%%%%%%%%%%%%%%%%%%%%%%%%%
% NEW SECTION %%%%%%%%%%%%%%%%%%%%%%%%%%%%%%%%%%%%%%%%%%%%%%%%%%%%%%%%
%%%%%%%%%%%%%%%%%%%%%%%%%%%%%%%%%%%%%%%%%%%%%%%%%%%%%%%%%%%%%%%%%%%%%%
\head\newsection Survey and introduction\endhead

Many of the well-known special functions, such as the Jacobi polynomials and
Bessel functions, satisfy addition formulas, which can be found in e.g.
\cite{\Aske}, \cite{\ErdeHTF},
\cite{\Wats}. Often there exists a group theoretic interpretation of
such an addition formula. This means that there exists a group $G$ and a
representation $t$ of $G$ in a Hilbert space $V$ such that for a suitable basis
$\{ e_n\}$ of $V$ the matrix elements $t_{n,m}\colon G\to\C$ defined by
$t_{n,m}(g)= \langle t(g)e_m,e_n\rangle$ are known in terms of special
functions. Then the homomorphism property
$$
t_{n,m}(gh) = \sum_p t_{n,p}(g) t_{p,m}(h),\qquad g,h\in G
\tag\eqname{\vglhomoprop}
$$
gives an addition formula for a suitably choosen basis of $V$ and for certain
elements $g,h\in G$. Usually we number the basis such that $t_{0,0}$ is left
and right invariant with respect to a certain subgroup $K$ (a spherical
function) and we use \thetag{\vglhomoprop} for $n=m=0$. See
Vilenkin \cite{\Vile} and Vilenkin and Klimyk \cite{\VileK} for more
information.

Before the advent of quantum groups and quantum algebras (or $q$-algebras)
only a few addition formulas were known for $q$-special functions. Addition
formulas for $q$-Krawtchouk polynomials and $q$-Hahn polynomials were
proved by Dunkl \cite{\Dunk} and Stanton \cite{\StanGD} in the late
70's using the interpretation of these polynomials on Chevalley groups over a 
finite field. See also \cite{\StanSIAM} for a partial result (a product formula).
This method is limited to $q$-Krawtchouk and $q$-Hahn
polynomials, cf. the list in \cite{\StanAKS}.
Also, by a purely analytic method Rahman and Verma \cite{\RahmV} have shown
that also the continuous $q$-ultraspherical polynomials satisfy an addition
formula similar to Gegenbauer's addition formula for the ultraspherical
polynomials.

With the introduction of quantum groups and quantum algebras (or quantised
universal enveloping algebras) a new setting for various $q$-special
functions emerged, see e.g. the survey papers by Koornwinder \cite{\KoorOPTA},
Noumi \cite{\Noum} and the author \cite{\KoelAAM} for the quantum group approach
and Floreanini and Vinet \cite{\FlorVLMP}, \cite{\FlorVAP} 
and references given later in this
section for the quantum algebra approach. For a good introduction to
quantum groups, quantised universal enveloping algebras and their applications
we refer to Chari and Pressley \cite{\CharP}, which also contains a
wealth of references on the subject.

From this setting two approaches for
deriving addition formulas for $q$-special functions can be obtained.
To understand the quantum group approach we recall that a quantum group
${\Cal A}_q(G)$ is a deformation (with deformation parameter $q$) of a function
algebra ${\Cal A}(G)$ on the
group $G$. For a compact group $G$ this would mean the polynomials on $G$ or,
in Woronowicz's \cite{\Woro} approach, the continuous functions on $G$.
We assume that the comultiplication
$$
\Delta \colon {\Cal A}(G)\to {\Cal A}(G)\otimes {\Cal A}(G)
\cong {\Cal A}(G\times G),
\qquad \bigl( \Delta f\bigr) (g,h) = f(gh),
$$
is well-defined and we demand that
$\Delta \colon {\Cal A}_q(G)\to {\Cal A}_q(G)\otimes {\Cal A}_q(G)$
survives the deformation unchanged (in a suitable sense). 
In the deformation the representations of $G$
correspond to corepresentations of ${\Cal A}_q(G)$ and for a suitable
basis of the representation space we get matrix elements
$t_{n,m}\in{\Cal A}_q(G)$ satisfying the cohomomorphism property
$$
\Delta(t_{n,m}) = \sum_p t_{n,p}\otimes t_{p,n},
\tag\eqname{\vglcohomoprop}
$$
which must be considered as the analogue of \thetag{\vglhomoprop}. If the
matrix elements can be expressed in terms of $q$-special functions we consider
\thetag{\vglcohomoprop} as an (implicit) addition formula. Usually
${\Cal A}_q(G)$ is a non-commutative algebra and thus we need to investigate
the representations of ${\Cal A}_q(G)$ in order to convert
\thetag{\vglcohomoprop} into an addition formula involving $q$-special
functions in commuting variables. Usually, $n=m=0$ in
\thetag{\vglcohomoprop} and $t_{0,0}$ can be considered as a spherical
function. In this paper we show that in the case of the quantum $SU(2)$ group 
we can obtain a general addition formula for Askey-Wilson polynomials from
\thetag{\vglcohomoprop} from which various known examples can be obtained as
special or limiting cases.

Let us now briefly consider the quantum algebra approach. A quantum algebra
is a deformation ${\Cal U}_q{\frak g}$ of the universal enveloping algebra
${\Cal U}{\frak g}$ of the Lie algebra ${\frak g}$. The representation theory
of ${\Cal U}_q{\frak g}$ is usually similar to the representation theory
of ${\frak g}$. Classically we can obtain elements of the corresponding
group $G$ by exponentiating Lie algebra elements. In the quantum algebra approach 
the action of $exp_q(\a_1X_1)\ldots exp_q(\a_nX_n)$ is calculated in a 
representation of ${\Cal U}_q{\frak g}$. Here $exp_q$ can be one of the 
$q$-analogues of the exponential function, $\a_i$ are scalars 
and $X_i$ are generators of ${\Cal U}_q{\frak g}$. For a 
suitable basis $\{ f_m\}$ of the representation space we get
$$
exp_q(\a_1X_1)\ldots exp_q(\a_nX_n)\, f_m = \sum_k U_{m,k}(\a_1,\ldots,\a_n)\, f_k.
\tag\eqname{\vglhomoqalg}
$$
The matrix coefficients $U_{m,k}(\a_1,\ldots,\a_n)$ can 
be calculated in terms of special functions in $\a_1,\ldots,\a_n$. By working with
explicit realisations, in which the $f_m$ correspond to certain special functions,
addition formulas can be derived from \thetag{\vglhomoqalg}.
This method is motivated by the classical relation between Lie algebras and special
functions as described in Miller's book \cite{\Mill} and references 
for the $q$-algebra approach are given later on.

It should be observed that in the quantum group case the coalgebra structure,
and in particular the comultiplication $\Delta$, is needed to find the addition
formula \thetag{\vglcohomoprop} in non-commuting variables. The algebra structure
of ${\Cal A}_q(G)$ is needed to transform \thetag{\vglcohomoprop} into an addition
formula for $q$-special functions in commuting variables, although identities for
$q$-special functions in non-commuting variables are of interest in their own right. 
In contrast, in the quantum
algebra approach only the algebra structure is needed, and there are examples of
quantum algebras used in relation with $q$-special functions which do not carry
a bialgebra structure.

Let us now give some references to the literature for several addition
formulas using one of these methods.
For addition formulas for $q$-Legendre polynomials
from the quantum group theoretic point of view see Koornwinder \cite{\KoorIM},
Masuda et. al \cite{\MasuMNNU}, Noumi and Mimachi \cite{\NoumMPJA,
\NoumMCMP, \NoumMDMJ, \NoumMCM}, Vaksman and Soibelman \cite{\VaksS}
for (implicit) addition formulas in non-commuting variables and
the author \cite{\KoelSIAM, \KoelCJM, \KoelCJMs}, Koornwinder \cite{\KoorSIAM}
(analytically proved by Rahman \cite{\RahmPAMS}) and
Noumi and Mimachi \cite{\NoumMPJA} for addition formulas involving only
commuting variables.
For addition formulas for $q$-Bessel functions from the
quantum group theoretic
point of view see the author \cite{\KoelDMJ, \KoelIM}. For the
$q$-algebra point of view see Floreanini and Vinet \cite{\FlorVJMP,
\FlorVCJP, \FlorVEsterel}, Kalnins and Miller \cite{\KalnM},
Kalnins, Miller and Mukherjee \cite{\KalnMiMuSIAM}. For analytic proofs
of related addition formulas for $q$-Bessel functions see
the author \cite{\KoelITSF, \KoelJCAM}, Koelink and Swarttouw \cite{\KoelS},
Koornwinder and Swarttouw \cite{\KoorS}, Rahman \cite{\RahmCJM},
Swarttouw \cite{\Swar}. Of these references \cite{\KoelITSF},
\cite{\KoorS}
are closely related to the $q$-algebra approach and
\cite{\KoelJCAM} is
related to the quantum group approach.
Using the quantum group approach it is possible to derive an
addition formula for $q$-disk polynomial involving little $q$-Jacobi polynomials,
cf. Floris \cite{\Flor}, Floris and Koelink \cite{\FlorK}. The announcement
\cite{\DijkN} by Dijkhuizen and Noumi suggests that generalisations to Askey-Wilson
polynomials might be possible in this context.
Using other quantum algebras it is possible to obtain addition
formulas for general ${}_r\vp_s$ (Floreanini and Vinet \cite{\FlorVJGTP}),
for basic Lauricella $\vp_D$ (Floreanini, Lapointe and Vinet \cite{\FlorLV}),
for continuous analogues of addition formulas (Floreanini and
Vinet \cite{\FlorVJCAM}, Kalnins and Miller \cite{\KalnM}),
for $q$-Laguerre polynomials (Kalnins, Manocha and Miller \cite{\KalnMM},
Kalnins and Miller \cite{\KalnM}, Kalnins, Miller and
Mukherjee \cite{\KalnMiMuJMP}). See Groza and Kachurik \cite{\GrozK}
for addition and product formulas from the quantum $SU(2)$ group interpretation
of $q$-Krawtchouk, $q$-Hahn and $q$-Racah polynomials.

In this paper we show how to derive an explicit addition formula for
Askey-Wilson polynomials involving $3$ parameters from the implicit addition
formula, i.e. involving non-commuting variables. So we start off with an
identity of the type as in \thetag{\vglcohomoprop} and we convert it into an
identity for Askey-Wilson polynomials. In \S 2 we recall results on the
quantum $SU(2)$ group and its relation with Askey-Wilson polynomials. In
\S 3 we consider suitable vectors in the representation space of an irreducible
representation of ${\Cal A}_q(SU(2))$ in which all the elements 
of ${\Cal A}_q(SU(2))$ under
consideration act as a multiplication operator or as a shift operator. 
From this we obtain
in \S 4 an addition formula for Askey-Wilson polynomials from which various
known addition formulas can be obtained as special and limit cases.

The notation for $q$-shifted factorials and $q$-hypergeometric series 
follows the excellent book \cite{\GaspR} by Gasper and Rahman.

\demo{Acknowledgement} I thank Mizan Rahman and Serge\u\i\ Suslov for sending
their preprint \cite{\AskeRS} and for answering a number of
questions. Most of the work for this paper was done at the Katholieke
Universiteit Leuven, and I thank Walter Van~Assche and Alfons Van~Daele
for their hospitality.\enddemo

%%%%%%%%%%%%%%%%%%%%%%%%%%%%%%%%%%%%%%%%%%%%%%%%%%%%%%%%%%%%%%%%%%%%%%
% NEW SECTION %%%%%%%%%%%%%%%%%%%%%%%%%%%%%%%%%%%%%%%%%%%%%%%%%%%%%%%%
%%%%%%%%%%%%%%%%%%%%%%%%%%%%%%%%%%%%%%%%%%%%%%%%%%%%%%%%%%%%%%%%%%%%%%
\head\newsection Generalised matrix elements on the quantum $SU(2)$
group\endhead

In this section we recall the relation between the quantum $SU(2)$ group
and the Askey-Wilson polynomials. We also give the appropriate version of 
\thetag{\vglcohomoprop} in this case, which is the starting point for the 
derivation of the addition formula. References for this section are 
\cite{\NoumMPJA, \NoumMLNM} and \cite{\KoelAAM}, from which the notation has been
taken and where further references can be found.

$\A$ is the complex unital associative algebra generated
by $\a$, $\b$, $\g$, $\d$ subject to the relations
$$
\gathered
\a\b =q\b\a ,\quad \a\g = q\g\a ,\quad \b\d = q\d\b ,\quad \g\d = q\d\g
,\\
\b\g =\g\b ,\quad \a\d -q\b\g = \d\a - q^{-1}\b\g =1
\endgathered
%\tag\eqname{\vglcommrelAq}
$$
for some constant $q\in\C$.
With a $\ast$-operator given by
$$
\a^\ast = \d , \quad \b^\ast = -q\g ,\quad \g^\ast = -q^{-1}\b ,\quad
\d^\ast = \a
%\tag\eqname{\vglsteropAq}
$$
the algebra $\A$ becomes a $\ast$-algebra for real $q$, and from now on
we fix $0<q<1$.

We use the following realisation of $\A$ by a non-faithful
representation.
An irreducible infinite dimensional $\ast$-representation
$\p$ of $\A$ in the Hilbert space
$\Hi$ with orthonormal basis $\{ e_n\mid n\in\Zp\}$ is given by
$$
\p(\a) e_n = \sqrt{1-q^{2n}} e_{n-1}, \quad
\p(\g )e_n = q^n e_n,
%\tag\eqname{\vgldefreprA}
$$
where we use the convention $e_{-p}=0$ for $p\in\N$. Note that
$-q\p(\g)=\p(\b)$. For all $\xi\in\A$, $\p(\xi)$ is a bounded operator
on $\Hi$. The one-dimensional $\ast$-representation
$\t_\phi\colon\A\to\C$ is defined by
$$
\t_\phi(\a)=e^{i\phi},\qquad \t_\phi(\g)=0.
\tag\eqname{\vglonedimrepA}
$$

The algebra $\A$ is an example of a Hopf $\ast$-algebra.
The comultiplication $\Delta$, which is a $\ast$-homomorphism of
$\A\to\A\otimes\A$, is given on the generators by
$$
\Delta(\a)=\a\otimes\a + \b\otimes\g ,\quad
\Delta(\g )=\g\otimes\a + \d\otimes\g .
\tag\eqname{\vgldefDeltaonAq}
$$

There exist elements $b^l_{i,j}(\t,\s)\in\A$, $l\in\hZp$, $i,j\in\{
-l,-l+1,\ldots,l-1,l\}$, $\s,\t\in\R\cup\{ \infty\}$ such that
$$
\Delta\bigl(b^l_{i,j}(\t,\s)\bigr) = \sum_{n=-l}^l
\bigl( D.b^l_{i,n}(\t,\m)\bigr) \otimes b^l_{n,j}(\m,\s),
\qquad\forall\;\m\in\R\cup\{\infty\}.
\tag\eqname{\vglabstractaddf}
$$
Here $D.\colon\A\to\A$ is the bijective algebra homomorphism given
by $D.\a=q^{-\hf}\a$, $D.\b=q^\hf\b$, $D.\g=q^{-\hf}\g$ and $D.\d=q^\hf\d$.
The case of most interest to us is $l\in\Zp$, $i=j=0$. So
\thetag{\vglabstractaddf} is an example of the cohomorphism property
\thetag{\vglcohomoprop}.

In order to describe the generalised matrix 
elements $b^l_{i,j}(\t,\s)$ in more detail we
introduce certain simple elements of $\A$.
In case $l=\hf$ the $2\times 2$-matrix $b^\hf (\t,\s)$ is given by
$$
{1\over{\sqrt{(1+q^{2\s})(1+q^{2\t})} }}
\pmatrix \a_{\t,\s},& \b_{\t,\s}\\ \g_{\t,\s},&\d_{\t,\s} \endpmatrix =
b^{1/2}(\t,\s) =
\pmatrix b^{1/2}_{-1/2,-1/2}, & b^{1/2}_{-1/2,1/2}\\
b^{1/2}_{1/2,-1/2}, & b^{1/2}_{1/2,1/2}\endpmatrix
$$
with
$$
\aligned
\a_{\t,\s} &= q^{1/2}\a-iq^{\s-1/2}\b +iq^{\t+1/2}\g + q^{\s+\t-1/2}\d, \\
\b_{\t,\s} &= -q^{\s+1/2}\a-iq^{-1/2}\b -iq^{\s+\t+1/2}\g + q^{\t-1/2}\d, \\
\g_{\t,\s} &= -q^{\t+1/2}\a+iq^{\t+\s-1/2}\b +iq^{1/2}\g + q^{\s-1/2}\d, \\
\d_{\t,\s} &= q^{\t+\s+1/2}\a+iq^{\t-1/2}\b -iq^{\s+1/2}\g + q^{-1/2}\d.
\endaligned
\tag\eqname{\vgldefabgdst}
$$
We let $\s,\t\in\R\cup\{\infty\}$, e.g.
$\a_{\infty,\s} = q^{1/2}\a-iq^{\s-1/2}\b$ and $\a_{\infty,\infty} =
q^{1/2}\a$. Then it follows that
$$
\gathered
\a_{\t,\s} = \a_{\t,\infty}+q^\s\b_{\t,\infty}, \quad
\b_{\t,\s} = \b_{\t,\infty}-q^\s\a_{\t,\infty}, \\
\g_{\t,\s} = \g_{\t,\infty}+q^\s\d_{\t,\infty}, \quad
\d_{\t,\s} = \d_{\t,\infty}-q^\s\g_{\t,\infty}.
\endgathered
\tag\eqname{\vglfactmateltlow}
$$
Next we define
$$
\aligned
\r &= {1\over 2}\bigr( q^{-\t-\s-1}\a_{\t+1,\s+1}\d_{\t,\s}
- q^{-\t-\s-1} - q^{\t+\s+1}\bigr) \\
&= {1\over 2}\bigr( q^{-\t-\s}\b_{\t+1,\s-1}\g_{\t,\s} +
q^{\s-\t-1} +q^{\t-\s+1}\bigr),
\endaligned
\tag\eqname{\vgldefrhosigmatau}
$$
cf. e.g. \cite{\KoelCJMs, prop.~3.3} for the last equality.
Note that $\r^\ast=\r$. The limit case $\s\to\infty$ of $\r$ is defined by
$\rti = \lim_{\s\to\infty} 2q^{\s+\t-1}\r
=q^{-1}\b_{\t,\infty}\g_{\t,\infty}+q^{2\t}$.

In order to be able to express $b_{i,j}^l(\t,\s)$ in these simple terms
we need Askey-Wilson polynomials $p_n(x;a,b,c,d|q)$, where we follow the
normalisation as in Askey and Wilson \cite{\AskeW, (1.15)}, see also
\cite{\GaspR, (7.5.2)}. We use the following notation for the Askey-Wilson
polynomials as $q$-analogues of the Jacobi polynomials;
$$
p_n^{(\a,\b)}(x;s,t|q) = p_n(x;q^{1/2}t/s,q^{1/2+\a}s/t,
-q^{1/2}/(st), -stq^{1/2+\b}|q).
\tag\eqname{\vgldefAWalsqJacobi}
$$
Observe that $\a,\b\to\infty$ in \thetag{\vgldefAWalsqJacobi} gives
Askey-Wilson polynomials with two parameters set to zero, which are the
Al-Salam--Chihara polynomials. So the Al-Salam--Chihara polynomials can
be considered as the Hermite case of the $q$-Jacobi polynomials
\thetag{\vgldefAWalsqJacobi}, and in this form they play a role
in the sequel.

The generalised matrix 
elements $b_{i,j}^l(\t,\s)$ can be expressed using these simple
elements and the Askey-Wilson polynomials.
Here we use the following four cases. For $n\in\Zp$ we have
$$
\aligned
b^l_{n,0}(\t,\s) &= d^{l,n}_{\t,\s}\, c^n_{n,0}(\t,\s)\,
p^{(n,n)}_{l-n}(\r;q^\t, q^\s | q^2), \\
b^l_{0,n}(\t,\s) &= d^{l,n}_{\s,\t}\, c^n_{0,n}(\t,\s)\,
p^{(n,n)}_{l-n}(\r;q^\s, q^\t | q^2),\\
b^l_{-n,0}(\t,\s) &= d^{l,n}_{-\t,-\s}\, c^{n}_{-n,0}(\t,\s)\,
p^{(n,n)}_{l-n}(\r;q^{-\t}, q^{-\s} | q^2), \\
b^l_{0,-n}(\t,\s) &= d^{l,n}_{-\s,-\t}\, c^{n}_{0,-n}(\t,\s)\,
p^{(n,n)}_{l-n}(\r;q^{-\s}, q^{-\t} | q^2)
\endaligned
%\tag\eqname{\vglblijts}
$$
with the so-called minimal elements in $\A$ defined by
$$
\aligned
c^n_{n,0}(\t,\s) &=  q^{n(1-\s)}
\prod_{j=0}^{n-1} \d_{\t+2n-1-2j,\s-1}\,\g_{\t+2n-2-2j,\s},\\
c^n_{0,n}(\t,\s) &=  q^{n(1-\t)}
\prod_{j=0}^{n-1} \d_{\t-1,\s+2n-1-2j}\,\b_{\t,\s+2n-2-2j},\\
c^n_{-n,0}(\t,\s) &=  q^{-n(\s+2\t-2n)}
\prod_{j=0}^{n-1} \b_{\t-2n+1+2j,\s-1}\,\a_{\t-2n+ 2+2j,\s},\\
c^n_{0,-n}(\t,\s) &=  q^{-n(2\s+\t-2n)}
\prod_{j=0}^{n-1} \g_{\t-1,\s+1-2n+2j}\,\a_{\t,\s+2-2n+2j},
\endaligned
\tag\eqname{\vgldefminelts}
$$
and with the constant given by
$$
\gather
d^{l,n}_{\t,\s} = {{C^{l,0}(\s)C^{l,n}(\t)q^{-l}}\over
{(q^{2l+2n+2};q^2)_{l-n}}} = d^{l,n}_{\t,-\s},
\tag\eqname{\vgldefdlnts}\\
C^{l,j}(\s) = q^{l+j} \left[ {{2l}\atop{l-j}}\right]_{q^2}^{1/2}
\left( {{1+q^{-4j-2\s}}\over{(1+q^{-2\s})
 (-q^{2-2\s};q^2)_{l-j} (-q^{2+2\s};q^2)_{l+j} }} \right)^{1/2}.
\endgather
$$
Note that $C^{l,j}(\s)=C^{l,-j}(-\s)$. Since we work in a non-commutative
algebra, we have to be careful about the ordering in the product. We use
the convention that $\prod_{i=0}^k \psi_i=\psi_0\psi_1\ldots\psi_k$.

The following identity in the $\A$ is
the starting point for the `addition formula'.
The identity is obtained from
\thetag{\vglabstractaddf} with $i=j=0$ after applying
$id\otimes\t_{\phi/2}$. We assume from now on that $\t,\s,\m\in\R$.
Explicitly,
$$
\aligned
&d^{l,0}_{\t,\s} p_l^{(0,0)}(\r^\phi;q^\t,q^\s|q^2) =
d^{l,0}_{\m,\t} d^{l,0}_{\m,\s} p_l^{(0,0)}(\cos\phi;q^\m,q^\s|q^2)
p_l^{(0,0)}(D.\rho_{\t,\m};q^\m,q^\s|q^2)\\
&+\sum_{n=1}^l d^{l,n}_{\m,\t} d^{l,n}_{\m,\s} e^{-in\phi}
(-q^{1+\s+\m}e^{i\phi},q^{1-\s+\m}e^{i\phi};q^2)_n
p_{l-n}^{(n,n)}(\cos\phi;q^\m,q^\s|q^2) \\
&\qquad\qquad\qquad \times \bigl(D.c^n_{0,n}(\t,\m)\bigr)
p_{l-n}^{(n,n)}(D.\rho_{\t,\m};q^\m,q^\t|q^2) \\
&+\sum_{n=1}^l d^{l,n}_{-\m,-\t} d^{l,n}_{-\m,-\s} e^{-in\phi}
(-q^{1-\s-\m}e^{i\phi},q^{1+\s-\m}e^{i\phi};q^2)_n
p_{l-n}^{(n,n)}(\cos\phi;q^{-\m},q^{-\s}|q^2)\\
&\qquad\qquad\qquad \times \bigl(D.c^n_{0,-n}(\t,\m)\bigr)
p_{l-n}^{(n,n)}(D.\rho_{\t,\m};q^{-\m},q^{-\t}|q^2),
\endaligned
\tag\eqname{\vglstartaddformone}
$$
where $\r^\phi=(id\otimes\t_{\phi/2})\Delta(\r)$. To this identity
in $\A$ we apply the $\ast$-re\-pre\-sen\-ta\-tion $\p$ in order to obtain an
identity in the space ${\Cal B}\bigl(\Hi\bigr)$ of bounded linear operators.

Note that if we apply the one-dimensional $\ast$-representation
$\t_{\theta/2}$ to \thetag{\vglstartaddformone} instead of $\p$, we
obtain the (degenerate) addition formula for
the Askey-Wilson polynomials, cf. Noumi and Mimachi
\cite{\NoumMPJA, thm.~4}, see also \cite{\KoelSIAM, (3.15)},
\cite{\KoelAAM, (8.1)}.

%%%%%%%%%%%%%%%%%%%%%%%%%%%%%%%%%%%%%%%%%%%%%%%%%%%%%%%%%%%%%%%%%%%%%%
% NEW SECTION %%%%%%%%%%%%%%%%%%%%%%%%%%%%%%%%%%%%%%%%%%%%%%%%%%%%%%%%
%%%%%%%%%%%%%%%%%%%%%%%%%%%%%%%%%%%%%%%%%%%%%%%%%%%%%%%%%%%%%%%%%%%%%%
\head\newsection Basis of the representation space \endhead

To turn \thetag{\vglstartaddformone} into an identity for $q$-special
functions we have to study the operators $\p(\r^\phi)$ and
$\p(D.\r)$. These operators are given as a five term recurrence operator,
but it turns out that $\Hi$ has an orthogonal basis of eigenvectors of 
$\p(\rti)$ in which these operators are tridiagonal. This gives the oppurtunity
to determine the eigenvectors in terms of orthogonal polynomials, and in this case
the Al-Salam--Chihara polynomials are involved. The action of the minimal elements
on these eigenvectors is given by certain shift operators, and thus the action of
each of the generalised matrix elements can be calculated.

So we first recall the following basis of $\Hi$ and
the action of certain operators in this basis, cf.
\cite{\KoelCJM, \S 4}, \cite{\KoelCJMs, \S 3}.

\proclaim{Proposition \theoremname{\propeigvectprti}}
{\rm (i)} $\Hi$ has an orthogonal basis of the form
$v_\l=v_\l(q^\t)$, where $\l=-q^{2n}$, $n\in\Zp$,
$\l=q^{2\t+2n}$, $n\in\Zp$. The vector $v_\l$ is an eigenvector of the
self-adjoint operator $\p(\rti)$ for the eigenvalue $\l$. Moreover,
with the normalisation $\langle v_\l,e_0\rangle = 1$ we have
$$
\align
\langle v_\l,v_\l\rangle &= q^{-2n} (q^2;q^2)_n
(-q^{2-2\t};q^2)_n (-q^{2\t};q^2)_\infty,\qquad \l=-q^{2n},\\
\langle v_\l,v_\l\rangle &= q^{-2n} (q^2;q^2)_n
(-q^{2+2\t};q^2)_n (-q^{-2\t};q^2)_\infty,\qquad \l=q^{2\t+2n}.
\endalign
$$
{\rm (ii)} For $\l=-q^{2n}$, $\l=q^{2\t+2n}$, $n\in\Zp$, we have
$$
\gather
\p(\a_{\t,\infty})v_\l(q^\t) = iq^{\hf-\t}(1+\l) v_{\l/q^2}(q^{\t-1}), \qquad
\p(\b_{\t,\infty})v_\l(q^\t) = iq^\hf v_\l (q^{\t-1}),\\
\p(\g_{\t,\infty})v_\l(q^\t) = iq^\hf (q^{2\t}-\l)v_{\l}(q^{\t+1}),\qquad
\p(\d_{\t,\infty})v_\l(q^\t) = -iq^{\hf+\t} v_{\l q^2}(q^{\t+1}),
\endgather
$$
with the convention $v_{-q^{-2}}(q^\t)=0=v_{q^{2\t-2}}(q^\t)$.
\endproclaim

\demo{Remark \theoremname{\remorthdecHilbert}} The basis of proposition
\thmref{\propeigvectprti} induces the orthogonal decomposition
$\Hi=V_1\oplus V_2$, where $V_1$, respectively $V_2$, is spanned by
$v_{-q^{2n}}$, $n\in\Zp$, respectively $v_{q^{2\t+2n}}$, $n\in\Zp$.
When needed we use $V_i^\t=V_i$ to stress the dependence on $\t$.
\enddemo

From \thetag{\vgldefDeltaonAq} and \thetag{\vglonedimrepA} we see that
$(id\otimes \t_{\phi/2})\circ \Delta$ multiplies $\a$ and $\g$ by
$e^{i\phi/2}$ and $\b$ and $\d$ by $e^{-i\phi/2}$. From this we see that
$(id\otimes \t_{\phi/2})\circ \Delta$ multiplies $\a_{\t,\infty}$ and
$\g_{\t,\infty}$ by $e^{i\phi/2}$ and $\b_{\t,\infty}$ and $\d_{\t,\infty}$
by $e^{-i\phi/2}$ and from this we find the action of
$(id\otimes \t_{\phi/2})\circ \Delta$ on $\a_{\t,\s}$, etcetera, cf.
\thetag{\vglfactmateltlow}.
Next recall \thetag{\vgldefrhosigmatau} to find
$$
\align
&\Bigl[ 2q^{\t+\s}(\p\otimes \t_{\phi/2})\Delta(\r)-q^{2\s-1}-q^{2\t+1}
\Bigr] v_\l =
(\p\otimes \t_{\phi/2})\Delta(\b_{\t+1,\s-1}\g_{\t,\s})v_\l \\
&\quad = \p\Bigl( (e^{-i\phi/2}\b_{\t+1,\infty}-
q^{\s-1}e^{i\phi/2}\a_{\t+1,\infty})
(e^{i\phi/2}\g_{\t,\infty}+ q^\s e^{-i\phi/2}\d_{\t,\infty})\Bigr)
v_\l,
\endalign
$$
and this can be calculated explicitly by proposition
\thmref{\propeigvectprti}(ii). We obtain
$$
\multline
2(\p\otimes\t_{\phi/2})\Delta(\r)v_\l = \\ e^{-i\phi}qv_{\l q^2} +
e^{i\phi} q^{-1} (1 - q^{-2\t}\l)(1+\l)v_{\l/q^2} +
\l q^{1-\t}(q^{-\s}-q^\s)v_\l.
\endmultline
\tag\eqname{\vglrecursiegeneraltwee}
$$

For convenience we restrict our attention from now on to the subspace
$V_1$ of $\Hi$. We use the notation $v_n=v_n(q^\t)$ for the orthogonal
basis $v_{-q^{2n}}=v_{-q^{2n}}(q^\t)$, $n\in\Zp$, of $V_1$. Put $\l=-q^{2n}$ in
\thetag{\vglrecursiegeneraltwee} and compare the recurrence relation
with the three-term recurrence relation for the Al-Salam--Chihara polynomials
$h_n(x)=h_n(x;s,t|q)=\lim_{\a,\b\to\infty} p_n^{(\a,\b)}(x;s,t|q)$,
cf. \cite{\AskeW, (1.24)}, \cite{\GaspR, (7.5.3)},
$$
2xh_n(x)=h_{n+1}(x) + (t-t^{-1})q^{n+\hf}s^{-1} h_n(x)
+ (1-q^n)(1+q^ns^{-2}) h_{n-1}(x),
\tag\eqname{\vglthreetermrecASC}
$$
to see that the operator $(\p\otimes\t_{\phi/2})\Delta(\r)$ can be realised
as the multiplication operator on a suitable weighted $L^2$-space. This is the
content of the following proposition, which follows from the spectral theory of
Jacobi matrices, cf. Berezanski\u\i\ \cite{\Bere, Ch.~VII, \S 1}, Dombrowski
\cite{\Domb}.

\proclaim{Proposition \theoremname{\propinterLtwospace}} Denote by $dm$ the
orthogonality measure for the Al-Salam--Chihara polynomials
$h_n(\cdot;q^\t,q^\s|q^2)$ normalised by
$\int_\R dm(x)=(-q^{2\t};q^2)_\infty$ and define the mapping
$\Lambda\colon V_1\to L^2(dm)$ by
$$
\Lambda\colon v_n(q^\t) \longmapsto
q^{-n} e^{in\phi}\, h_n(\cdot;q^\t,q^\s|q^2).
$$
Then $\Lambda$ extends to a unitary operator and
$\Lambda\bigl((\p\otimes\t_{\phi/2})\Delta(\r)\bigr) = M\Lambda$,
where the multiplication operator $M\colon L^2(dm)\to L^2(dm)$ is
defined by $Mf\colon x\mapsto xf(x)$.
\endproclaim

\demo{Remark \theoremname{\rempropunitary}} The fact that $\Lambda$ extends
to a unitary operator follows from the determinacy of the moment problem for
the Al-Salam--Chihara polynomials. The multiplication
operator $M$ is a bounded operator on $L^2(dm)$, since the support of the
orthogonality measure $dm$ is compact.
\enddemo

Next we study $\p(D.\r)$. In a completely analogous way we prove
$$
2\p(D.\r)v_\l = q^2v_{\l q^2} +
q^{-2} (1 - q^{-2\t}\l)(1+\l)v_{\l/q^2} +
\l q^{1-\t}(q^{-\s}-q^\s)v_\l.
%\tag\eqname{\vglrecursiegeneraldrie}
$$
Restricting to $V_1$ we see that we can find eigenvectors of $\p(D.\r)$ for
the eigenvalue $y$ of the form $\sum_{n=0}^\infty p_n(y)v_n$ if and only if
the $p_n$'s satisfy
$$
\align
2yp_n(y) &=q^{-2}(1+q^{2n+2-2\t})(1-q^{2n+2})p_{n+1}(y)  \\
&\qquad\qquad\qquad  + q^{1-\t+2n}(q^\s-q^{-\s})p_n(y) + q^2 p_{n-1}(y), \\
2yp_0(y) &=q^{-2}(1+q^{2-2\t})(1-q^2)p_1(y) +
q^{1-\t}(q^\s-q^{-\s})p_0(y),
\endalign
$$
so that the eigenvector is completely determined by $p_0(y)$.
Here we use the convention that $v_{-1}=0$.
So we find the eigenvector
$$
u_y=u_y(q^\t,q^\s)=\sum_{n=0}^\infty {{q^{2n}
h_n(y;q^\t,q^\s|q^2)}\over{(q^2,-q^{2-2\t};q^2)_n}}\,
v_n(q^\t), \qquad y\in{\text{supp}}(dm),
\tag\eqname{\vgleigvectDrts}
$$
by using the three-term recurrence relation \thetag{\vglthreetermrecASC}.
To see that $u_y\in\Hi$ for $y$ in the support of $dm$
we can use the fact that the similarly defined vector
with $q^{2n}$ replaced by $q^n$ is a generalised eigenvector of the
self-adjoint operator $\p(\r)$ and then use \cite{\Bere, Ch.~VII, \S 1.1,
(1.24)}. Another way to see this is to use the asymptotic properties of the
Al-Salam--Chihara polynomials as $n\to\infty$, cf. e.g. Askey and Ismail
\cite{\AskeI, \S 3.1}.

\proclaim{Proposition \theoremname{\propactioDabgdonux}} The operators
$\p(D.\a_{\t,\s})$ and $\p(D.\b_{\t,\s})$ map $V_1^\t$ into $V_1^{\t-1}$
and the operators $\p(D.\g_{\t,\s})$ and $\p(D.\d_{\t,\s})$ map $V_1^\t$
into $V_1^{\t+1}$. Moreover,
$$
\align
\p(D.\a_{\t,\s})\, u_y(q^\t,q^\s) &=
{{iq^{1+\s}}\over{1+q^{2-2\t}}} (1+2yq^{1-\t-\s}+q^{2-2\t-2\s})
\,u_y(q^{\t-1},q^{\s-1}), \\
\p(D.\b_{\t,\s})\, u_y(q^\t,q^\s) &=
{{iq}\over{1+q^{2-2\t}}} (1-2yq^{1+\s-\t}+q^{2+2\s-2\t})
\,u_y(q^{\t-1},q^{\s+1}),\\
\p(D.\g_{\t,\s})\, u_y(q^\t,q^\s) &= i(1+q^{2\t})
\,u_y(q^{\t+1},q^{\s-1}), \\
\p(D.\d_{\t,\s})\, u_y(q^\t,q^\s) &= -iq^\s(1+q^{2\t})
\,u_y(q^{\t+1},q^{\s+1}).
\endalign
$$
\endproclaim

\demo{Proof} First observe that $D.\g_{\t,\infty}=q^{-\hf}\g_{\t,\infty}$
and $D.\d_{\t,\infty}=q^\hf\d_{\t,\infty}$. Now \thetag{\vgldefabgdst}
and proposition \thmref{\propeigvectprti}(ii) show that
$\p(D.\g_{\t,\s})$ maps $V_1^\t$ into $V_1^{\t+1}$. And similarly for the
other operators.

To prove the second statement we note that the results for $D.\g_{\t,\s}$ and
$D.\d_{\t,\s}$ imply the results for $D.\b_{\t,\s}$ and $D.\a_{\t,\s}$ by
\thetag{\vgldefrhosigmatau}, since $u_y$ is an eigenvector of
$\p(D.\r)$. Recall, cf. \cite{\KoelCJMs, cor.~3.3},
$$
(D.\g_{\t,\s})(D.\r)=(D.\rho_{\t+1,\s-1})(D.\g_{\t,\s}).
$$
Since the eigenspaces of $\p(D.\r)$ in $V_1^\t$
are one-dimensional, this implies
$$
\p(D.\g_{\t,\s})u_y(q^\t,q^\s) = C\, u_y(q^{\t+1},q^{\s-1})
$$
for some constant $C$. To calculate $C$ we take the inner product with
$v_0(q^{\t+1})$ and use that $\p$ is a $\ast$-representation to get
$$
\align
&C\, (-q^{2+2\t};q^2)_\infty = \langle u_y(q^\t,q^\s),
\p\bigl( (D.\g_{\t,\s})^\ast\bigr) v_0(q^{\t+1})\rangle \\
&\quad = \langle u_y(q^\t,q^\s),
\p\bigl(-q^{-\hf}\b_{\t+1,\infty}+q^{\s-\hf}\a_{\t+1,\infty}\bigr)
v_0(q^{\t+1})\rangle = i (-q^{2\t};q^2)_\infty
\endalign
$$
by \thetag{\vglfactmateltlow}, $\g_{\t,\infty}^\ast=-\b_{\t+1,\infty}$,
$\d_{\t,\infty}^\ast=q^{-1}\a_{\t+1,\infty}$ and proposition
\thmref{\propeigvectprti}(ii). This implies the value for $C$. The
statement for $\p(D.\d_{\t,\s})$ is proved analogously.
\qed\enddemo

\demo{Remark \theoremname{\remidforAWpols}} We can rewrite the result of
proposition \thmref{\propactioDabgdonux} as identities for Al-Salam--Chihara
polynomials by use of \thetag{\vgleigvectDrts},
\thetag{\vglfactmateltlow} and proposition \thmref{\propeigvectprti}(ii) to
find special cases of
$$
\multline
(1-q^{2n}abcd)(1-2ay+a^2)p_n(y;aq,b,c,d| q) = \\-ap_{n+1}(y;a,b,c,d| q)
+(1-abq^n)(1-acq^n)(1-adq^n)p_n(y;a,b,c,d| q),
\endmultline
$$
which can be proved from the orthogonality relations for the Askey-Wilson
polynomials, and of
$$
\multline
(1-q^{2n-2}abcd)p_n(y;a/q,b,c,d|q) = (1-q^{n-2}abcd)p_n(y;a,b,c,d|q)\\
-{a\over q}(1-q^n)(1-bcq^{n-1})(1-bdq^{n-1})(1-cdq^{n-1})p_{n-1}(y;a,b,c,d|q)
\endmultline
$$
which is a special case of the connection coefficients in
\cite{\AskeW, \S 6}.
\enddemo

From proposition \thmref{\propactioDabgdonux} and \thetag{\vgldefminelts}
we obtain the following corollary by iteration.

\proclaim{Corollary \theoremname{\coractionminelts}} With $y=(z+z^{-1})/2$
we have
$$
\align
\p\bigl(D.c_{0,n}^n(\t,\s)\bigr)\, u_y(q^\t,q^\s) &=
q^{n(\t+\s+n)} (q^{1+\s-\t}z, q^{1+\s-\t}/z;q^2)_n\, u_y(q^\t,q^{\s+2n}), \\
\p\bigl(D.c_{0,-n}^n(\t,\s)\bigr)\, u_y(q^\t,q^\s) &= (-1)^n
q^{n(\t-\s+n)} \\
&\qquad\qquad \times (-q^{1-\s-\t}z, -q^{1-\s-\t}/z;q^2)_n
\,u_y(q^\t,q^{\s-2n}).
\endalign
$$
\endproclaim

%%%%%%%%%%%%%%%%%%%%%%%%%%%%%%%%%%%%%%%%%%%%%%%%%%%%%%%%%%%%%%%%%%%%%%
% NEW SECTION %%%%%%%%%%%%%%%%%%%%%%%%%%%%%%%%%%%%%%%%%%%%%%%%%%%%%%%%
%%%%%%%%%%%%%%%%%%%%%%%%%%%%%%%%%%%%%%%%%%%%%%%%%%%%%%%%%%%%%%%%%%%%%%
\head\newsection Addition formula for Askey-Wilson polynomials \endhead

The results of section 3 allow us to calculate the action of the elements in
\thetag{\vglstartaddformone} under the representation $\p$ on suitable vectors.
In this section we show how we can obtain a very general type of addition formula
for the Legendre case, i.e. $\a=\b=0$, of the $q$-Jacobi polynomials defined 
in \thetag{\vgldefAWalsqJacobi}. We also show how this formula covers known addition
formulas for $q$-Legendre polynomials and we end with some open problems.

Apply $\p$ to \thetag{\vglstartaddformone} and let the resulting
identity in ${\Cal B}\bigl(\Hi\bigr)$ act on the eigenvector
$u_y(q^\t,q^\m)$ of $\p(D.\rho_{\t,\m})$. The action of the right hand side
of \thetag{\vglstartaddformone} follows from corollary
\thmref{\coractionminelts}. To the resulting identity we apply the
unitary operator $\Lambda$ of proposition \thmref{\propinterLtwospace},
which also shows how the left hand side of \thetag{\vglstartaddformone} looks
under $\Lambda$.
In order to find the result we need the action of $\Lambda$ on the right hand side of 
\thetag{\vglstartaddformone}. For this we have
to calculate the following function of
$x$ in the Hilbert space $L^2(dm)$;
$$
\bigl(\Lambda u_y(q^\t,q^\m)\bigr)(x) =
\sum_{m=0}^\infty {{q^m e^{im\phi}}\over{(q^2,-q^{2-2\t};q^2)_m}}
h_m(y;q^\t,q^\m|q^2) h_m(x;q^\t,q^\s|q^2).
\tag\eqname{\vglinnerprodeigvect}
$$
Note that in case $\m=\s$ this is the Poisson kernel for the Al-Salam--Chihara
polynomials at $t=qe^{i\phi}$.
This absolutely convergent expression, which can be considered as a
non-symmetric Poisson kernel for the Al-Salam--Chihara polynomials,
has been evaluated by Askey, Rahman and
Suslov \cite{\AskeRS, \S 14, Case~II} in terms of a very-well-poised
${}_8\vp_7$-series. Explicitly, if we define for $|t|<1$, $x=\cos\psi$,
$y=\cos\theta$ the function $P(t;x,y;\t;\s,\m) =$
$$
\multline
{{(-q^{-\s-\m}t, q^{1+\s-\t}te^{i\theta},
q^{1+\s-\t}te^{-i\theta}, q^{1+\m-\t}te^{i\psi},
q^{1+\m-\t}te^{-i\psi};q^2)_\infty}\over {(q^{2+\m+\s-2\t}t,
te^{i\theta+i\psi}, te^{i\theta-i\psi}, te^{i\psi-i\theta},
te^{-i\theta-i\psi};q^2)_\infty}} \\
\times {}_8W_7 (q^{\m+\s-2\t}t;-q^{\s+\m}t, q^{1+\m-\t}e^{i\theta},
q^{1+\m-\t}e^{-i\theta}, \\
q^{1+\s-\t}e^{i\psi}, q^{1+\s-\t}e^{-i\psi}; q^2 ,-q^{-\s-\m}t),
\endmultline
\tag\eqname{\vgl8W7Poisson}
$$
then $\bigl(\Lambda u_y(q^\t,q^\m)\bigr)(x) = P(qe^{i\phi};x,y;\t;\s,\m)$.
Here we use the standard notation for very-well-poised basic hypergeometric
series as in \cite{\GaspR, \S 2.1};
$$
{}_8W_7(a;b,c,d,e,f;q,z) = {}_8\vp_7\left(
{{a\, ,q\sqrt{a}\, ,-q\sqrt{a}\, ,b\, ,c\, ,d\, ,e\, ,f}\atop
{\sqrt{a},-\sqrt{a},qa/b,qa/c,qa/d,qa/e,qa/f}};q,z\right).
$$
These remarks prove the following addition theorem for the Askey-Wilson
polynomials.

\proclaim{Theorem \theoremname{\thmaddformAWgen}}
We have the
following `addition formula' for Askey-Wilson polynomials;
$$
\aligned
&d^{l,0}_{\t,\s}
p_l^{(0,0)}(x;q^\t,q^\s|q^2)
P(qe^{i\phi};x,y;\t;\s,\m) = \\ &\qquad
A_0\, p_l^{(0,0)}(z;q^\m,q^\s|q^2)
p_l^{(0,0)}(y;q^\m,q^\t|q^2)
P(qe^{i\phi};x,y;\t;\s,\m)\\
&+\sum_{n=1}^l  A_n\, e^{-in\phi}
(-q^{1+\s+\m}e^{i\phi},q^{1-\s+\m}e^{i\phi};q^2)_n
(q^{1+\m-\t}e^{i\theta}, q^{1+\m-\t}e^{-i\theta};q^2)_n \\
&\qquad\times p_{l-n}^{(n,n)}(z;q^\m,q^\s|q^2)
p_{l-n}^{(n,n)}(y;q^\m,q^\t|q^2)
P(qe^{i\phi};x,y;\t;\s,\m+2n)\\
&+\sum_{n=1}^l B_n\, e^{-in\phi}
(-q^{1-\s-\m}e^{i\phi},q^{1+\s-\m}e^{i\phi};q^2)_n
(-q^{1-\m-\t}e^{i\theta}, -q^{1-\m-\t}e^{-i\theta};q^2)_n \\
&\qquad\times p_{l-n}^{(n,n)}(z;q^{-\m},q^{-\s}|q^2)
p_{l-n}^{(n,n)}(y;q^{-\m},q^{-\t}|q^2)
P(qe^{i\phi};x,y;\t;\s,\m-2n)
\endaligned
$$
with $A_n=d^{l,n}_{\m,\t} d^{l,n}_{\m,\s} q^{n(\t+\m+n)}$,
$B_n= (-1)^nd^{l,n}_{-\m,-\t} d^{l,n}_{-\m,-\s} q^{n(\t-\m+n)}$,
$x=\cos\psi$, $y=\cos\theta$, $z=\cos\phi$ and $d^{l,n}_{\m,\s}$ and
the Poisson kernel defined by \thetag{\vgldefdlnts} and
\thetag{\vgl8W7Poisson}.
\endproclaim

Initially theorem \thmref{\thmaddformAWgen}
only holds for $x=\cos\psi$ almost everywhere as an identity in
$L^2(dm)$, but by continuity it holds everywhere.

\demo{Remark \theoremname{\remotherspace}} Theorem
\thmref{\thmaddformAWgen} has been derived from the identity
\thetag{\vglstartaddformone} after applying $\p$ to it. We have only
considered the resulting operator identity restricted to the subspace
$V_1^\t$ of $\Hi$. If we study the restriction to the subspace $V_2^\t$, we
obtain the same formula after replacing $\t$, $\s$ and $\m$ by
$-\t$, $-\s$ and $-\m$. Replacing $\p$ by any other irreducible
infinite-dimensional representation of $\A$ does not lead to greater
generality.
\enddemo

\demo{Remark \theoremname{\remspecialcases}} (i) To recover the
(degenerate) addition formula for the 2-pa\-ra\-me\-ter family of Askey-Wilson
polynomials, i.e. the result of applying $\t_{\theta/2}$ to
\thetag{\vglstartaddformone}, we formally replace $e^{i\theta}$ by
$q^{-1}e^{i\theta}$, $e^{i\phi}$ by $q^{-1}e^{i\phi}$, and
$e^{i\psi}$ by $q^{-1}e^{i\theta+i\phi}$. For this choice of parameters
the very-well-poised ${}_8\vp_7$-series in the non-symmetric Poisson
kernel reduces to a very-well-poised ${}_6\vp_5$-series, which can be
summed by \cite{\GaspR, (II.20)}.

(ii) In \cite{\KoelCJMs} the case $i=j=0$, $\m\to\infty$ of
\thetag{\vglabstractaddf} is turned into an addition formula for the
2-parameter family of Askey-Wilson polynomials involving big $q$-Jacobi
polynomials in the sum and the role of the very-well-poised
${}_8\vp_7$-series is taken over by $q$-Laguerre polynomials, i.e. the
polynomials obtained from \thetag{\vgldefAWalsqJacobi} by letting
$\b\to\infty$. In \cite{\KoelCJMs} $\p\otimes\p$ is used as the
representation of $\A\otimes\A$ instead of $\p\otimes\t_{\phi/2}$ as
here. It is possible to obtain this result formally as a limit case of
theorem \thmref{\thmaddformAWgen} as follows. Replace $e^{i\phi}$ by
$q^{\m+\s-1}z^{-1}$ and let $\m\to\infty$ such that
$2q^{\m+\t-1}\cos\theta$ tends to $y$. Observe that then
$2q^{\m+\s-1}\cos\phi$ tends to $z$. Then we can use the limit
transition of the Askey-Wilson polynomials to the big $q$-Jacobi
polynomials as described by Koornwinder \cite{\KoorZSE, prop.~6.1}.
It remains to consider the limit case $\m\to\infty$ of the non-symmetric
Poisson kernel \thetag{\vgl8W7Poisson}, which can be done directly in
the very-well-poised ${}_8\vp_7$-series. The result is a
${}_3\vp_2$-series which can be written as a $q$-Laguerre polynomial
when $y$ and $z$ are mass points of the orthogonality measure for the
corresponding big $q$-Jacobi polynomials by simple transformations for
${}_3\vp_2$-series.

The limit case can also be evaluated directly from the infinite sum
\thetag{\vglinnerprodeigvect}. Using the limit as described in
\cite{\KoorZSE, prop.~6.1} we get
$$
q^{m\m}h_m(\cos\theta;q^\t,q^{\m+2n}|q^2)\to y^mq^{(1-\t)m}
(-y^{-1}q^{-2n};q^2)_m
$$
as $\m\to\infty$. The resulting sum is then a generating function for the
Al-Salam--Chihara polynomials, which has been evaluated by Serge\u\i\
Suslov (private communication). Explicitly, for the Al-Salam--Chihara
polynomial defined by
$$
s_m(\cos\psi;a,b,|q) = a^{-m} (ab;q)_m \ {}_3\vp_2\left(
{{q^{-m},ae^{i\psi},ae^{-i\psi}}\atop{ab,\ 0}};q,q\right)
$$
we have the generating function
$$
\sum_{m=0}^\infty {{(u;q)_m t^m}\over{(q,ab;q)_m}} s_m(\cos\psi;a,b|q)
= {{(ute^{-i\psi};q)_\infty}\over{(te^{-i\psi};q)_\infty}}
\ {}_3\vp_2\left( {{u,ae^{i\psi},ae^{-i\psi}}\atop{ab,ute^{-i\psi}}};q,
te^{i\psi}\right).
$$
Again for the choices of $y$ and $z$ as above this can be rewritten as
a $q$-Laguerre polynomial.
\enddemo

In general it seems hard to reduce theorem \thmref{\thmaddformAWgen} to a
polynomial identity, but for the special case $\s=\t=\m=0$ we can obtain
the following addition formula for the continuous $q$-Legendre
polynomials, cf. Rahman and Verma \cite{\RahmV, (1.24)}, see also
\cite{\KoelSIAM, thm.~4.1}.

\proclaim{Corollary \theoremname{\coraddformcontLeg}} For the
continuous $q$-ultraspherical polynomials defined by
$$
C_n(\cos\psi;\b|q) = \sum_{k=0}^n {{(\b;q)_k(\b;q)_{n-k}}\over
{(q;q)_k(q;q)_{n-k}}}e^{i(n-2k)\psi}
$$
we have
$$
\multline
C_l(\cos\psi;q^2\vert q^4) = q^l C_l(\cos\phi;q^2\vert q^4)
C_l(\cos\theta;q^2\vert q^4) + \\
\sum_{n=1}^l  {{q^{l-2n}(1+q^{4n})(q^2;q^4)_n}\over
{(-q^2;q^2)_{2n}(q^4;q^4)_n}}
\left[{{l+n}\atop{2n}}\right]_{q^4}^{-1}
C_{l-n}(\cos\theta;q^{2+4n}\vert q^4)C_{l-n}(\cos\phi;q^{2+4n}\vert q^4)\\
\times
e^{-in(\theta+\phi)} (q^2 e^{2i\theta}, q^2 e^{2i\phi};q^4)_n
\, {}_4\varphi_3 \left( {{q^{-4n},q^{4n},qe^{i(\phi+\theta)} e^{i\psi},
qe^{i(\phi+\theta)} e^{-i\psi}}\atop{ q^2,\; q^2e^{2i\phi},\;
q^2e^{2i\theta}}};q^4,q^4\right).
\endmultline
$$
\endproclaim

\demo{Proof}
We can apply the same reductions as in \cite{\KoelSIAM, \S 4} to theorem
\thmref{\thmaddformAWgen}. First observe
$$
P(t;\cos\theta,\cos\psi;0;0,0) = {{(t^2;q^4)_\infty}\over{
(te^{i\psi+i\theta}, te^{i\psi-i\theta}, te^{i\theta-i\psi},
te^{-i\psi-i\theta};q^4)_\infty}}, \quad |t|<1.
$$
The last identity can be observed from $h_n(x;1,1|q^2)=
H_n(x|q^4)$, (use e.g. the three-term recurrence relation) and
using Rogers's expression for the Poisson kernel of the continuous
$q$-Hermite polynomials \cite{\AskeIqH}, \cite{\Bres} or by applying
the summation formula \cite{\GaspR, (II.18)}. Now the corollary follows
from the following lemma, which has been proved by Mizan Rahman
(private communication) using
transformation and summation theorems for basic hypergeometric series.

\proclaim{Lemma \theoremname{\lemredforcontqLeg}} For $n\in\Zp$ and the
Poisson kernel $P$ defined by \thetag{\vgl8W7Poisson} we have \par
$q^{n^2}(qe^{i\theta},q^{-i\theta};q^2)_n
P(qe^{i\phi};\cos\psi,\cos\theta;0;0,2n)$
$$\multline
+(-1)^n q^{n^2}(-qe^{i\theta},-q^{-i\theta};q^2)_n
P(qe^{i\phi};\cos\psi,\cos\theta;0;0,-2n) = \\
{{2e^{-in\theta}(q^2e^{i\theta};q^4)_n(q^2e^{2i\phi};q^4)_\infty}\over{
(qe^{i\phi+i\psi+i\theta}, qe^{i\phi+i\psi-i\theta},
qe^{i\phi+i\theta-i\psi}, qe^{i\phi-i\psi-i\theta};q^4)_\infty}} \\
\times\,
{}_4\varphi_3 \left( {{q^{-4n},q^{4n},qe^{i(\theta+\phi)} e^{i\psi},
qe^{i(\theta+\phi)} e^{-i\psi}}\atop{ q^2,\; q^2e^{2i\phi},\;
q^2e^{2i\theta}}};q^4,q^4\right).
\endmultline
$$
\endproclaim

Of course, it also works the other way round; assuming theorem
\thmref{\thmaddformAWgen} and corollary \thmref{\coraddformcontLeg}
gives lemma \thmref{\lemredforcontqLeg}. Since Rahman's proof is limited
to the case $\s=\t=\m=0$ we do not give it here. \qed\enddemo

\demo{Remark \theoremname{\remopenproblems}} (i) Because of the similarity
between \thetag{\vglabstractaddf} and the group theoretic proof of the
addition formula for the Legendre polynomials, cf.
\cite{\Vile, \S III.4.2}, \cite{\VileK, \S 6.6.1}, we want to have a limit
transition of theorem \thmref{\thmaddformAWgen} to the addition formula
for the Legendre polynomials as $q\uparrow 1$. A straighforward limit
does not seem possible in general, but it does work for the Rahman-Verma addition
formula in corollary \thmref{\coraddformcontLeg}, cf. \cite{\RahmV}.
Also the more sophisticated technique of Van
Assche and Koornwinder \cite{\VAsscK}, which can be used to handle the result for the
case $\m=\infty$, is not applicable, since we do not
have a three-term recurrence relation for the non-symmetric Poisson
kernel in $\m$, i.e. a simple relation between the non-symmetric Poisson
kernels with $\m-2$, $\m$ and $\m+2$. \par
(ii) Usually an addition formula leads to a product formula. In
this case this can be done if $P(qe^{i\phi};x,y;\t;\s,\m+2n)$ is
part of a set of e.g. biorthogonal rational functions in $x$ with
respect to $n$. This seems not known.
\enddemo

%%%%%%%%%%%%%%%%%%%%%%%%%%%%%%%%%%%%%%%%%%%%%%%%%%%%%%%%%%%%%%%%%%%%%%%%%
%%%%%%%%%%%%%%%%%References%%%%%%%%%%%%%%%%%%%%%%%%%%%%%%%%%%%%%%%%%%%%%%
%%%%%%%%%%%%%%%%%%%%%%%%%%%%%%%%%%%%%%%%%%%%%%%%%%%%%%%%%%%%%%%%%%%%%%%%%


\Refs

\ref\no\Aske
\by Askey, R.
\yr 1975
\book Orthogonal Polynomials and Special Functions
\bookinfo CBMS-NSF Regional Conference Series Applied Math. {\bf 21}
\publaddr SIAM, Philadelphia PA
\endref

\ref\no\AskeIqH
\by Askey, R., and Ismail, M.E.H.
\yr 1983
\paper A generalization of ultraspherical polynomials
\inbook Studies in Pure Mathematics
\ed Erd\H{o}s, P.
\publaddr Birkh\"auser, Basel
\pages 55--78
\endref

\ref\no\AskeI
\bysame
\yr 1984
\paper Recurrence relations, continued fractions and orthogonal
polynomials
\jour Memoirs Amer. Math. Soc.
\vol 49
\issue 300
\endref

\ref\no\AskeRS
\by Askey, R.A., Rahman, M., and Suslov, S.K.
\yr 1994
\paper On a general $q$-Fourier transformation with nonsymmetric kernels
\paperinfo preprint, series~2, no.~21, Carleton University
\endref

\ref\no\AskeW
\by Askey, R., and Wilson, J.
\yr 1985
\paper Some basic hypergeometric orthogonal polynomials that
generalize Jacobi polynomials
\jour Memoirs Amer. Math. Soc.
\vol 54
\issue 319
\endref

\ref\no\Bere
\by Berezanski\u\i, J.M.
\yr 1968
\book Expansions in Eigenfunctions of Selfadjoint Operators
\bookinfo Transl. Math. Monographs 17
\publaddr Amer. Math. Soc., Providenc RI
\endref

\ref\no\Bres
\by Bressoud, D.M.
\yr 1980
\paper A simple proof of Mehler's formula for $q$-Hermite polynomials
\jour Indiana Univ. Math. J.
\vol 29
\pages 577--580
\endref

\ref\no\CharP
\by Chari, V., and Pressley, A.
\yr 1994
\book A Guide to Quantum Groups
\publaddr Cambridge University Press, Cambridge
\endref

\ref\no\DijkN
\yr 1995
\by Dijkhuizen, M.S., and Noumi, M.
\paper A family of quantum projective spaces and related $q$-hypergemeotric
orthogonal polynomials
\paperinfo announcement
\endref

\ref\no\Domb
\by Dombrowski, J.
\yr 1990
\paper Orthogonal polynomials and functional analysis
\inbook Orthogonal Polynomials: Theory and Practice
\ed Nevai, P.
\bookinfo NATO ASI series C, vol. 294
\publaddr Kluwer, Dordrecht
\pages 147--161
\endref

\ref\no\Dunk
\by Dunkl, C.F.
\yr 1977
\paper An addition theorem for some $q$-Hahn polynomials
\jour Monatsh. Math.
\vol 85
\pages 5--37
\endref

\ref\no\ErdeHTF
\by Erd\'elyi,A., Magnus, W., Oberhettinger, F., and Tricomi, F.G.
\yr 1953, 1955
\book Higher Transcendental Functions
\bookinfo 3 volumes
\publaddr McGraw-Hill
\endref

\ref\no\FlorLV
\by Floreanini, R., Lapointe, L., and Vinet, L.
\yr 1994
\paper A quantum algebra approach to basic mulivariable special functions
\jour J. Phys. A: Math. Gen.
\vol 27
\pages 6781--6797
\endref

\ref\no\FlorVJMP
\by Floreanini, R., and Vinet, L.
\yr 1992
\paper Addition formulas for $q$-Bessel functions
\jour J. Math. Phys.
\vol 33
\pages 2984--2988
\endref

\ref\no\FlorVPLA
\bysame
\yr 1992
\paper Using quantum algebras in $q$-special function theory
\jour Phys. Lett. A
\vol 170
\pages 21--28
\endref

\ref\no\FlorVLMP
\bysame
\yr 1993
\paper On the quantum group and quantum algebra approach to $q$-special
functions
\jour Lett. Math. Phys.
\vol 27
\pages 179--190
\endref

\ref\no\FlorVJGTP
\bysame
\yr 1993
\paper An algebraic interpretation of the $q$-hypergeometric function
\jour J. Group Theory Phys.
\vol 1
\pages 1--10
\endref

\ref\no\FlorVAP
\bysame
\yr 1993
\paper Quantum algebras and $q$-special functions
\jour Ann. Phys.
\vol 221
\pages 53--70
\endref

\ref\no\FlorVCJP
\bysame
\yr 1994
\paper Generalized $q$-Bessel functions
\jour Canad. J. Phys.
\vol 72
\pages 345--354
\endref

\ref\no\FlorVJCAM
\bysame
\yr 1995
\paper $q$-Gamma and $q$-beta functions in quantum algebra representation
theory
\jour J. Comp. Appl. Math.
\toappear
\endref

\ref\no\FlorVEsterel
\bysame
\yr 1995
\paper Basic Bessel functions and $q$-difference equations
\inbook Proc. `Symmetries and integrability properties of difference
equations'
\toappear
\endref

\ref\no\Flor
\by Floris, P.G.A.
\yr 1994
\paper Addition formula for $q$-disk polynomials
\paperinfo report W-94-24, University of Leiden
\endref

\ref\no\FlorK
\by Floris, P.G.A., and Koelink, H.T.
\yr 1995
\paper Addition formula for little $q$-disk polynomials
\paperinfo in preparation
\endref

\ref\no\GaspR
\by Gasper, G., and Rahman, M.
\yr 1990
\book Basic Hypergeometric Series
\bookinfo Encyclopedia Math. Appl. 35
\publaddr Cambridge University Press, Cambridge
\endref

\ref\no\GrozK
\by Groza, V.A., and Kachurik, I.I.
\yr 1990
\paper Addition and multiplication theorems for Krawtchouk, Hahn and
Racah $q$-polynomials
\paperinfo (in Russian)
\jour Dokl. Akad. Nauk Ukrain SSR, Ser. A
\vol 89
\pages 3--6
\endref

\ref\no\KalnMM
\by Kalnins, E.G., Manocha, H.L., and Miller, W.
\yr 1992
\paper Models of $q$-algebra representations: tensor products of special
unitary and oscillator algebras
\jour J. Math. Phys.
\vol 33
\pages 2365--2383
\endref

\ref\no\KalnM
\by Kalnins, E.G., and Miller, W.
\yr 1994
\paper Models of $q$-algebra representations: $q$-integral transforms
and addition theorems
\jour J. Math. Phys.
\vol 35
\pages 1951--1975
\endref

\ref\no\KalnMiMuJMP
\by Kalnins, E.G., Miller, W., and Mukherjee, S.
\yr 1993
\paper Models of $q$-algebra representations: matrix elements of the
$q$-oscillator algebra
\jour J. Math. Phys.
\vol 34
\pages 5333--5356
\endref

\ref\no\KalnMiMuSIAM
\bysame
\yr 1994
\paper Models of $q$-algebra representations: the group of plane motions
\jour SIAM J. Math. Anal.
\vol 25
\pages 513--527
\endref

\ref\no\KoelITSF
\by Koelink, H.T.
\yr 1993
\paper A basic analogue of Graf's addition formula and related
formulas
\jour Integral Transforms and Special Functions
\vol 1
\pages 165--182
\endref

\ref\no\KoelSIAM
\bysame
\yr 1994
\paper The addition formula for
continuous $q$-Legendre polynomials and associated spherical
elements on the $SU(2)$ quantum group related to Askey-Wilson
polynomials
\jour SIAM J. Math. Anal.
\vol 25
\pages 197--217
\endref

\ref\no\KoelDMJ
\bysame
\yr 1994
\paper The quantum group of plane motions and the Hahn-Exton $q$-Bessel
function
\jour Duke Math. J.
\vol 76
\pages 483--508
\endref

\ref\no\KoelCJM
\bysame
\yr 1995
\paper Addition formula for
big $q$-Legendre polynomials from the quantum $SU(2)$ group
\jour Canad. J. Math.
\toappear
\endref

\ref\no\KoelIM
\bysame
\yr 1995
\paper The quantum group of plane motions and basic Bessel
functions
\jour Indag. Math.
\toappear
\endref

\ref\no\KoelAAM
\bysame
\yr 1995
\paper Askey-Wilson polynomials and the quantum $SU(2)$ group: survey
and applications
\jour Acta Appl. Math.
\toappear
\endref

\ref\no\KoelJCAM
\bysame
\yr 1995
\paper Yet another basic analogue of Graf's addition formula
\jour J. Comp. Appl. Math.
\toappear
\endref

\ref\no\KoelCJMs
\bysame
\yr 1994
\paper Addition formula for $2$-parameter family of Askey-Wilson polynomials
\paperinfo preprint K.U. Leuven
\endref

\ref\no\KoelS
\by Koelink, H.T., and Swarttouw, R.F.
\yr 1995
\paper A $q$-analogue of Graf's addition formula for the Hahn-Exton
$q$-Bessel function
\jour J. Approx. Theory
\vol 81
\pages 260--273
\endref

\ref\no\KoorIM
\by Koornwinder, T.H.
\yr 1989
\paper Representations of the twisted $SU(2)$ quantum group and some
$q$-hypergeometric orthogonal polynomials
\jour Proc. Kon. Ned. Akad. van Wetensch., Ser.~A {\bf 92} (Indag. Math.
{\bf 51})
\pages 97--117
\endref

\ref\no\KoorOPTA
\bysame
\yr 1990
\paper Orthogonal polynomials in connection with quantum groups
\inbook Orthogonal Polynomials: Theory and Practice
\ed Nevai, P.
\bookinfo NATO ASI series C, vol. 294
\publaddr Kluwer, Dordrecht
\pages 257--292
\endref

\ref\no\KoorSIAM
\bysame
\yr 1991
\paper The addition formula for
little $q$-Legendre polynomials and the $SU(2)$ quantum group
\jour SIAM J. Math. Anal.
\vol 22
\pages 295--301
\endref

\ref\no\KoorZSE
\bysame
\yr 1993
\paper Askey-Wilson polynomials
as zonal spherical functions on the $SU(2)$ quantum group
\jour SIAM J. Math. Anal.
\vol 24
\pages 795--813
\endref

\ref\no\KoorS
\by Koornwinder, T.H., and Swarttouw, R.F.
\yr 1992
\paper On $q$-analogues of the Fourier and Hankel transforms
\jour Trans. Amer. Math. Soc.
\vol 333
\pages 445--461
\endref

\ref\no\Mill
\by Miller, W.
\yr 1968
\book Lie Theory and Special Functions
\publaddr Academic Press, New York
\endref

\ref\no\MasuMNNU
\by Masuda, T., Mimachi, K., Nakagami, Y., Noumi, M., and Ueno, K.
\yr 1991
\paper Representations of the quantum group $SU_q(2)$ and the little
$q$-Jacobi polynomials
\jour J. Funct. Anal.
\vol 99
\pages 357--386
\endref

\ref\no\Noum
\by Noumi, M.
\yr 1991
\paper Quantum groups and $q$-orthogonal polynomials. Towards a
realization of Askey-Wilson polynomials on $SU_q(2)$
\inbook Special Functions
\eds Kashiwara, M., and Miwa, T.
\bookinfo ICM-90 Satellite Conference Proceedings
\publaddr Springer-Verlag, New York
\endref

\ref\no\NoumMPJA
\by Noumi, M., and Mimachi, K.
\yr 1990
\paper Askey-Wilson polynomials and the quantum group $SU_q(2)$
\jour Proc. Japan Acad., Ser. A
\vol 66
\pages 146--149
\endref

\ref\no\NoumMCMP
\bysame
\yr 1990
\paper Quantum $2$-spheres and big $q$-Jacobi polynomials
\jour Comm. Math. Phys.
\vol 128
\pages 521--531
\endref

\ref\no\NoumMDMJ
\bysame
\yr 1991
\paper Rogers's $q$-ultraspherical polynomials on a quantum $2$-sphere
\jour Duke Math. J.
\vol 63
\pages 65--80
\endref

\ref\no\NoumMCM
\bysame
\yr 1992
\paper Spherical functions on a family of quantum $3$-spheres
\jour Comp. Math.
\vol 83
\pages 19--42
\endref

\ref\no\NoumMLNM
\bysame
\yr 1992
\paper Askey-Wilson polynomials as spherical functions on $SU_q(2)$
\inbook Quantum Groups
\bookinfo Lecture Notes Math. 1510
\ed Kulish, P.P.
\publaddr Springer-Verlag, New York
\pages 98--103
\endref

\ref\no\RahmCJM
\by Rahman, M.
\yr 1988
\paper An addition theorem and some product formulas for $q$-Bessel
functions
\jour Canad. J. Math.
\vol 40
\pages 1203--1221
\endref

\ref\no\RahmPAMS
\bysame
\yr 1989
\paper A simple proof of Koornwinder's addition formula for the little
$q$-Legendre polynomials
\jour Proc. Amer. Math. Soc.
\vol 107
\pages 373--381
\endref

\ref\no\RahmV
\by Rahman, M., and Verma, A.
\yr 1986
\paper Product and addition formulas for the continuous
$q$-ultra\-sphe\-ri\-cal polynomials
\jour SIAM J. Math. Anal.
\vol 17
\pages 1461--1474
\endref

\ref\no\StanSIAM
\by Stanton, D.
\yr 1980
\paper Product formulas for $q$-Hahn polynomials
\jour SIAM J. Math. Anal.
\vol 11
\pages 100--107
\endref

\ref\no\StanGD
\bysame
\yr 1981
\paper Three addition theorems for some $q$-Krawtchouk polynomials
\jour Geom. Dedicata
\vol 10
\pages 403--425
\endref

\ref\no\StanAKS
\bysame
\yr 1984
\paper Orthogonal polynomials and Chevalley groups
\inbook Special Functions: Group Theoretical Aspects and Applications
\ed Askey, R.A., Koornwinder, T.H., and Schempp, W.
\publaddr Reidel, Dordrecht
\pages 87--128
\endref

\ref\no\Swar
\by Swarttouw, R.F.
\yr 1992
\paper An addition theorem and some product formulas for the Hahn-Exton
$q$-Bessel functions
\jour Canad. J. Math.
\vol 44
\pages 867--879
\endref

\ref\no\VaksS
\by Vaksman, L.L., and Soibelman, Y.S.
\yr 1988
\paper Algebra of functions on the quantum group $SU(2)$
\jour Funct. Anal. Appl.
\vol 22
\pages 170--181
\endref

\ref\no\VAsscK
\by Van Assche, W., and Koornwinder, T.H.
\yr 1991
\paper Asymptotic behaviour for Wall polynomials and the addition formula
for little $q$-Legendre polynomials
\jour SIAM J. Math. Anal.
\vol 22
\pages 302--311
\endref

\ref\no\Vile
\by Vilenkin, N.J.
\yr 1968
\book Special Functions and the Theory of Group Representations
\bookinfo Transl. Math. Monographs 22
\publaddr Amer. Math. Soc., Providenc RI
\endref

\ref\no\VileK
\by Vilenkin, N.J., and Klimyk, A.U.
\yr 1991, 1993
\book Representation of Lie Groups and Special Functions
\bookinfo 3 volumes
\publaddr Kluwer, Dordrecht
\endref

\ref\no\Wats
\by Watson, G.N.
\yr 1944
\book A Treatise on the Theory of Bessel Functions
\bookinfo 2nd ed.
\publaddr Cambridge University Press, Cambridgr
\endref

\ref\no\Woro
\by Woronowicz, S.L.
\yr 1987
\paper Compact matrix pseudogroups
\jour Comm. Math. Phys.
\vol 111
\pages 613--665
\endref

\endRefs
\enddocument